 \documentclass[11pt]{article}
 \usepackage[top=1in,bottom=1in,left=1.5in,right=1.5in]{geometry}
  \usepackage{amsmath,amssymb}
  \usepackage{latexsym}% defines $\Box$ for LaTeX2e
 \usepackage[dvips]{pstricks} % PSTricks
  \usepackage{pst-node}
  \usepackage[dvips]{graphicx}
  \newcommand{\C}{\mathbb{C}}

  \newcommand{\R}{\mathbb{R}}
  \newcommand{\Z}{\mathbb{Z}}

  \newcommand{\x}{\mathbf{x}}

  \newcommand{\0}{\mathbf{0}}
  \newcommand{\1}{\mathbf{1}}

  \newcommand{\cM}{\mathcal{M}}

  \newcommand{\hs}{\hspace*{\parindent}}
  \newcommand{\proof}{\hs \textbf{Proof.\ }}

  \newcommand{\trans}{^\top}
  \newcommand{\qed}{\hspace*{\fill} $\Box$\\}
  
  \newcommand{\haff}{\mathrm{haf}}

  \newcommand{\rS}{\mathrm{S}}

  \newcommand{\perm}{\mathrm{perm}}

  \newtheorem{theo}{\bfseries \hs Theorem}[section]

  \newtheorem{lemma}[theo]{\bfseries \hs Lemma}

  \numberwithin{equation}{section} % Automatically number equations within sections

\begin{document}

 \title{Analogs of the van der Waerden and Tverberg\\conjectures for haffnians}

 \author
 {Shmuel Friedland\\
 Department of Mathematics, Statistics and Computer Science\\
 University of Illinois at Chicago\\
 Chicago, Illinois 60607-7045, USA\\
 \emph{email}:friedlan@uic.edu}

 \date{February 17, 2011}
 \maketitle

 \begin{abstract}
 We discuss here analogs of van der Waerden and Tverberg permanent conjectures for haffnians on the convex set
 of matrices whose extreme points are symmetric permutation matrices with zero diagonal.

 \end{abstract}



 \section{Introduction}  Let $G=(V,E)$ be an undirected graph with no loops but possibly with multi-edges.
 Denote by $\cM_k(G)$ the subset of all $k$-matchings.  It is well known that for a given graph $G$ there exist polynomial time algorithms to find if $\cM_k(G)$ is empty or not \cite{Ed65}.
 However, even for bipartite graph the problem of finding $\#\cM_k(G)$ is $\#$P-complete problem, see \cite{Val} for perfect matching
 and \cite{FL6} for $k$-matchings.  There are fully randomized polynomial approximation schemes to estimate $\#\cM_k(G)$ in bipartite graphs \cite{JSV,FL6}.

 The aim of this paper is to discuss lower bounds for $\#\cM_k(G)$
 for regular graphs.  A lot is known for regular bipartite graphs.  This follows from the validity of the van der Waerden and Tverberg conjectures for permanents of doubly stochastic matrices \cite{vdW,Tve,Fr79,Ego,Fal,Fr82}.  Improved lower bounds were obtained for regular bipartite graphs \cite{Vor, Sch,Gur1,FG08}.
 One of the main features of regular bipartite graphs that they can be expressed as a $k$ edge-disjoint union of perfect matches.  More precisely the Birkhoff
 theorem yields that any doubly stochastic matrix is a convex combination of permutation matrices.

 For regular nonbipartite graph the situation is more complex.  There exist simple cubic graphs which do not have a perfect matching.
 See for example the Sylvester graph \cite[Figure P.3, p'xii]{LP}.
 The celebrated Petersen's theorem claims that a cubic graph with two bridges at most has a perfect matching.  It was only shown very recently that
 a simple cubic graph without a bridge have exponential number of perfect matchings in the number of vertices \cite{EKKKN}.
 Let $G=(V,E)$ be an undirected graph with no loops and denote by $A(G)$ the adjacency matrix of $G$.  Then $A(G)$ is a symmetric matrix with zero diagonal
 and nonnegative integer off-diagonal entries.  For a set $T\subset \R$ denote by $\rS_0(n,T)$ the set of symmetric matrices with zero diagonal whose off-diagonal entries are in $T$.  Thus, $B=[b_{ij}]\in\rS_0(n,\R_+)$ can be viewed as viewed as a weighted complete graph $K_n$, where $0\le b_{ij}$ is the weight of the edge $(i,j)$.  Assume that $n=2m$ is even.  Denote by $\haff(B)$, the \emph{haffnian} of $B$, the sum of the weighted perfected matches of $K_n$, given by $B$.  For any $n$ and positive integer $k$, such that $2k\le n$, and $B\in \rS_0(n,\R_+)$, denote by $\haff_k(B)$ the sum of  all $2k$ haffnians of principle
 submatrices of order $2k$ of $B$.  Equivalently, $\haff_k(B)$ is the sum of all weighted $k$-matches in $K_n$ given by $B$.
 Thus for a given graph $G=(V,E)$, $\haff_k(A(G))$ is the number of $k$-matches in $G$.

 Suppose that $G=(V,G)$ is bipartite, where $V=V_1\cup V_2, E\subset V_1\times V_2$.  Then $\hat A(G)\in \Z_+^{|V_1|\times |V_2|}$ is the bipartite adjacency
 matrix of $G$, i.e. $A(G)=\left[\begin{array}{cc}0&\hat A(G)\\(\hat A(G))\trans& 0\end{array}\right]$.  More generally, given positive integers $l,m$ and an $l\times m$ nonnegative matrix $C=[c_{ij}]\in \R_+^{l\times m}$, then $C$ can be viewed as a weighted bipartite graph $K_{l,m}$.  Thus $B=\left[\begin{array}{cc}0&C\\C\trans& 0\end{array}\right]$ be the adjacency matrix of the weighted graph $K_{l,m}$.  Then $\haff_k(B)=\perm_k(C)$, where
 $\perm_k(C)$ is the sum of permanents of all $k\times k$ submatrices of $C$.  Denote by $\Omega_n\subset \R_+^{n\times n}$ the set of doubly stochastic matrices.  The proved Tverberg conjecture states \cite{Fr82}.
 \begin{equation}\label{tvercon}
 \min_{C\in \Omega_n} \perm_k(C)= \perm_k(\frac{1}{n}\hat A(K_{n,n})) \textrm{ for } k=2,\ldots,n.
 \end{equation}
 Equality  holds if and only $C=\frac{1}{n}\hat A(K_{n,n})$.  The case $k=n$ is the van der Waerden conjecture.
 (\ref{tvercon}) for $k=n$ implies immediately that any $r$-regular bipartite graph on $2n$ vertices has ar least $(\frac{r}{e})^n$ perfect matchings.
 For better bounds see \cite{Vor,Sch,Gur1}.

 The main success in proving the Tverberg conjecture and its sharper analogs for $r$-regular bipartite graphs
 can be attributed to the notion of \emph{hyperbolicity}.  Suffices to say that a product of linear factors is hyperbolic.
 Hence for any nonnegative matrix $C\in \R_+^{n\times n}$ the polynomial $f(\x):=\prod_{i=1}^n (C\x)_i$ is positive hyperbolic.
 ($f(\x)$ is a sum of monomial with nonnegative coefficients.)   Furthermore $\perm (C)=\frac{\partial^n}{\partial x_1\ldots\partial x_n} f$,
 is the mixed derivative of $f$.

 The aim of this paper is to introduce analogous problems to the van der Waerden and Tverberg conjectures.  First, one needs to introduce an analog notion to the notion
 of doubly stochastic matrices $\Psi_{2n}\subset \rS(2n,\R_+)\cap\Omega_{2n}$.  Namely, this is the convex set of symmetric doubly stochastic matrices with
 zero diagonal, whose extreme points are symmetric permutation matrices with zero diagonal.  $\Psi_{2n}$ was characterized by Edmonds \cite{Ed65a}, see \cite[Theorem 6.3, 2nd Proof, page 209]{CCPS} for a simple proof.
 Namely, it is the set of all $(2n)\times (2n)$ stochastic matrices $B=[b_{ij}]\in\Omega_{2n}$, which are symmetric and have zero diagonal, that satisfy the condition
 \begin{equation}\label{edmcon}
 \sum_{i,j\in S} b_{ij}\le |S|-1, \textrm{ for each } S\subset{1,\ldots,2n}, |S| \textrm{ odd and } 3\le |S|\le 2n-3.
 \end{equation}

 Our problem is to find or give a good lower bound for 
 \begin{equation}\label{fridprob}
 \min_{B\in\Psi_{2n}}\haff_k(B)= \mu_{k,n} \textrm{ for } k=2,\ldots,n.
 \end{equation}
 It is tempting to state, as in the case of the van der Waerden and Tverberg's conjectures that
 \begin{equation}\label{fridcon}
 \mu_{k,n}=\haff_k(\frac{1}{2n-1}A(K_{2n}))\textrm{ for } k=2,\ldots,n.
 \end{equation}
 Equality holds if and only if $B=\frac{1}{2n-1}A(K_{2n})$.
 See \cite{Fr11}.  
 (According to a recent e-mail from Leonid Gurvits, he stated this conjecture for $k=n$ in correspondence with E. Lieb on September 21, 2005.)
 As we show in \S3, this conjecture is true for $k=2$.  However for $k=n$ and $n$ big enough \eqref{fridcon} is wrong as explained below.

 Note that if $G=(V,E)$ is a $r$-regular graph without loops on an even number of vertices, then $\frac{1}{r}A(G)$ is in $\Psi_{|V|}$
 if and only if any vertex cut $S\subset V$ with an odd number of vertices, $3\le |S|\le |V|-3$, has at least $r$ edges.  Hence if Conjecture \eqref{fridcon} holds, then such a regular graph has at least $(\frac{r}{e})^{\frac{|V|}{2}}$ perfect matchings, see \eqref{approx}. 
 In \cite{CPS} the authors construct an infinite family of 3-edge connected graph $G=(V,E)$, i.e. an edge disjoint union of $3$-perfect matchings, for which the number of matchings is less then $c_F|V|(\frac{1+\sqrt{5}}{2})^{\frac{|V|}{12}}$. (Here $|V|=12k+4$ and $k=1,2,\ldots$.)  As $(\frac{1+\sqrt{5}}{2})^{\frac{1}{12}}< 1.017< \sqrt{\frac{3}{e}}\approx 1.05$ we must have that $\mu_{n,n}<\haff (\frac{1}{2n-1}A(K_{2n}))$ for $n\gg 1$.  (I would like to thank S. Norin for pointing out to me this fact.)
 
 Since $\mu_{n,n}$ is the minimum of the haffnian function, it follows that $\mu_{n+m,n+m}\le \mu_{n,n}\mu_{m,m}$.  Hence the sequence $\log\mu_{n,n}$ is subadditive. In particular the following limit exists
 \begin{equation}\label{defmu}
 \mu:=\lim_{n\to\infty} \frac{\log\mu_{n,n}}{n}.
 \end{equation}
  A weak analog of the van der Waerden conjecture is the claim that $\mu>-\infty$.  Note that the above example in \cite{CPS} implies that $\mu\le \frac{\log\frac{1+\sqrt{5}}{2}}{6}-\log 3$.  Other generalizations of the van der Waerden conjectures for perfect matchings in hypergraphs are
 considered in \cite{BS}.

 In \S2 we estimate evaluate $\haff_k(\frac{1}{2n-1}A(K_{2n}))$ and estimate the value of $\haff_n(\frac{1}{2n-1}A(K_{2n}))$
 for large $n$.  Using the notion of hyperbolicity we give a good lower bound of $\haff_k(B)$ for $B\in\Psi_{2n}$ with one positive
 eigenvalue.  In \S3 we show that $\frac{1}{2n-1}A(K_{2n})$ is a strict local minimum of $\haff_k(\cdot)$ on $\Psi_{2n}$ for $k=2,\ldots,n$.

 \section{Some equalities and lower estimates}
 \begin{lemma}\label{exactfor}  For positive integers $2\le k\le n$ we have
 \begin{eqnarray}\label{exactfor1h}
 \haff_k(\frac{1}{2n-1} A(K_{2n}))=\frac{1}{(2n-1)^k} {2n\choose 2k} \frac{1}{k!}\prod_{j=0}^{k-1} {2k-2j\choose 2},\\
 \haff_k(\frac{1}{n} A(K_{n,n}))=\frac{1}{n^k} {n \choose k}^2 k!\label{exactfor1p}.
 \end{eqnarray}
 In particular
 \begin{equation}\label{exactfor2}
  e^{-n}\sqrt{2}<\haff_n(\frac{1}{2n-1} A(K_{2n}))=\frac{(2n)!}{(2n-1)^n 2^n n!}<\haff_n(\frac{1}{n} A(K_{n,n}))=\frac{n!}{n^n}
 \end{equation}
 For $n\gg 1$ we have the following approximations
 \begin{equation}\label{approx}
 \haff_n(\frac{1}{2n-1} A(K_{2n}))\approx e^{-n}\sqrt{2e} , \quad \haff_n(\frac{1}{n} A(K_{n,n}))\approx e^{-n}\sqrt{2\pi n}.
 \end{equation}
 \end{lemma}
 \proof  We first compute the number of perfect matchings in $K_{2n}$.  For $l=1,\ldots, n$ choose the $l-th$ match between two distinct vertices
 of $K_{2n}$ and remove these two vertices from the vertices of $K_{2n}$.  Then number of choices for the $l-th$ pair is $2k-2(l-1)\choose 2$.
 Hence the total number of choices of $n$ pairs, taking in account the order of choices is $\prod_{j=0}^{n-1} {2n-2j\choose 2}$.
 Hence
 \[\haff_n(A(K_{2n}))=\frac{\prod_{j=0}^{n-1} {2n-2j\choose 2}}{n!}=\frac{(2n)!}{2^n n!}.\]
 Observe next that to find all $k$ matchings in $K_{2n}$ we first choose $2k$ vertices out of $2n$ vertices in $K_{2n}$, to obtain a subgraph $K_{2k}$ of $K_{2n}$.
 Then we compute all perfect matchings in $K_{2k}$.  So $\haff_k(A(K_{2n}))={2n\choose 2k} \haff_k(A(K_{2k}))$.
 This proves \eqref{exactfor1h}.  \eqref{exactfor1p} is well known \cite{Fr82}, and proved similarly.

 To obtain the first inequality in \eqref{exactfor2} we need the exact form of Stirling's formula \cite[p.52]{Fel}.
 \begin{equation}\label{stirling}
 m!=\sqrt{2\pi m}\; m^m e^{-m} e^{\frac{\theta_m}{12 m}}, \quad \theta_m \in (0,1), \quad m=1,2,\ldots .
 \end{equation}
 Hence
 \[\frac{(2n)!}{(2n-1)^n 2^n n!}\ge \frac{\sqrt{2\pi 2n}(2n)^{2n}e^n}{e^{2n} (2n-1)^n 2^n \sqrt{2\pi n}n^n e^{\frac{1}{12}}}= e^{-n} \sqrt{2}(\frac{2n}{2n-1})^n e^{-\frac{1}{12}}.\]
 Recall that the sequence $(\frac{m}{m-1})^{m-1}$ is an increasing sequence for $m=2,\ldots,$.
 Hence $(\frac{2n}{2n-1})^{n-\frac{1}{2}}$ is an increasing sequence.  As $n\ge 2$ we deduce that the left hand-side of the above inequality
 is greater than $e^{-n}\sqrt{2}(\frac{4}{3})^{\frac{3}{2}} e^{-\frac{1}{12}}> e^{-n}\sqrt{2}$.  This establishes the left-hand side of \eqref{exactfor2}.  To establish the second inequality of \eqref{exactfor2} divide the middle expression of \eqref{exactfor2} by its right-hand side to obtain $(\frac{2n}{2n-1})^n 4^{-n} {2n \choose n}$.  Recall that the sequence $(\frac{m}{m-1})^m$ is decreasing for $m=2,\ldots,$.
 Hence the sequence $(\frac{2n}{2n-1})^n$ is a decreasing sequence.  Hence to show the second inequality of \eqref{exactfor2}
 one needs to show that $(\frac{4}{3})^2 4^{-n} {2n \choose n}< 1$ for $n\ge 2$.  This claim follows easily by induction.
 \eqref{approx} follows straightforward from Stirling's formula.  \qed

 \begin{theo}\label{hypest}  Let $B\in \Psi_{2n}$ . Assume that $B\in\Psi_{2n}$ has exactly one positive eigenvalue.  Then
 \begin{equation}\label{hypest1}
 \haff_n(B)\ge (\frac{n-1}{n})^{(n-1)n}\approx e^{-n}\sqrt{e}.
 \end{equation}
 Moreover, for each $k=2,\ldots,n-1$
 \begin{equation}\label{hypest2}
 \haff_k(B)\ge \frac{(2n)^{2n-2k} (2n-k)! (2n)^k}{(2n-2k)!(2n-k)^{2n-k}2^k k!} (\frac{(2n-k-1)}{2n-k})^{(2n-k-1)k}.
 \end{equation}

 \end{theo}
 \proof The above inequalities follows from the results of \cite{FG06} as follows.  It is well known that the quadratic polynomial
 $\x\trans B\x$ is hyperbolic, if and only if it has exactly one positive eigenvalue, e.g. \cite[Lemma 6.1]{FG06}.
 Observe next that
 \begin{equation}\label{capb}
 \textrm{Cap} ((\x\trans B\x)^k):=\inf_{\x=(x_1,\ldots,x_{2n})\trans >\0} \frac{(\x\trans B \x)^k}{(\prod_{i=1}^{2n} x_i)^{\frac{k}{n}}}=(\1\trans B\1)^k=(2n)^k
 \end{equation}
 for each $B\in\Psi_{2n}$.
 Here $\1=(1,\ldots,1)\trans\in \R^{2n}$.  In view of arithmetic-geometric inequality we deduce that above equality if $B$ corresponds to a symmetric permutation matrix with zero diagonal.  Hence the above equality holds each $B\in\Psi_{2n}$.  Recall that
 \begin{equation}\label{haffkfor}
 \haff_k C=(2^{k} k!)^{-1} \sum_{1\le i_1 <\ldots<i_{2k}\le 2n}
 \frac{\partial ^{2k} }{\partial
 x_{i_1}\ldots\partial x_{i_{2k}}} (\x\trans C\x)^k, \; C\in\rS_{m}(0,\R_+),
 \end{equation}
 where $2k\le m$.  See \cite[\S6]{FG06}.
 Apply \cite[Theorem 3.1]{FG06} to the hyperbolic polynomial $\frac{(\x\trans B\x)^n}{2^n n!}$ to estimate its mixed derivative with respect to all
 $2n$ variables to deduce the inequality in \eqref{hypest1}.  (Replace $n$ by $2n$, note that $r_1=\ldots=r_{2n}=n$, $k=n+1$ and use the equality \eqref{capb}.)  Use the approximate expansion $\log(1-x)=-x -\frac{x^2}{2} +O(x^3)$ to deduce the approximation in \eqref{hypest1} for $n\gg 1$.
 \eqref{hypest2} is deduced similarly from \cite[Theorem 3.1]{FG06}.  \qed

 \section{Local conditions}\label{sec:loccon}
 \begin{lemma}\label{locmin}  Let $2\le k\le n$ be integers.  Then $C=\frac{1}{2n-1}A(K_{2n})$ is an interior point of
 the convex set $\Psi_{2n}$.  $C$ is a critical point of $\haff_k(\cdot)$
 on $\Psi_{2n}$.  The hessian of $\haff_k(\cdot)$ at $C$ has positive eigenvalues.  Hence $C$ is a unique local minimum of
 $\haff_k(\cdot)$.  In particular, Conjecture \ref{fridcon} holds for $k=2$.
 \end{lemma}
 \proof  Let $S\subset\{1,\ldots,2n\}$.  Then the sum of all the entries of $C$ for $i,j\in S$ is $\frac{|S|(|S|-1)}{2n-1}$, which is strictly less than $|S|-1$ if $|S|\in[2,2n-2]$.  Hence $C$ is an interior point $\Psi_{2n}$.  More precisely, each point $X\in\Psi_{2n}$ in the neighborhood of $C$ if and only if it is of the form $C+Y$, where $Y=[y_{ij}]$ belongs to the subspace
 \begin{equation}\label{defPhi}
 \Phi_{2n}=\{Y\in \rS_0(2n,\R),\; Y\1=\0\}.
 \end{equation}
 Observe next that for each integer $k\in [2,n]$ we have the equality
 \begin{equation}\label{haffkexp}
 \haff_k(F)=\frac{1}{k}\sum_{i< j} f_{ij}\haff_{k-1} (F[[2n]\setminus\{i,j\}]), \quad F=[f_{ij}]\in\rS_0(2n,\R).
 \end{equation}
 Here for any positive integer $l$ we let $[l]=\{1,\ldots,l\}$.  Furthermore, for any $X=[x_{ij}]\in\R^{l\times l}$ and any $S\subseteq [l]$
 we denote by $X[S]$ the principle submatrix $[x_{ij}]_{i,j\in S}$.

 Hence
 \begin{eqnarray*}
 \haff_k(C+Y)=\haff_k(C)+\frac{1}{2k}\sum_{i\ne j}y_{ij}\frac{1}{(2n-1)^{k-1}}\haff_{k-1}(A(K_{2n-2}))+\\
 \frac{2}{k(k-1)}\sum_{1\le i<j<p<q\le 2n\}|=4} (y_{ij}y_{pq}+y_{ip}y_{jq}+y_{iq}y_{jp}) \frac{1}{(2n-1)^{k-2}}h_{k-2}(A(K_{2n-4}))+O(\|X\|^3).
 \end{eqnarray*}
 Here $\haff_0(F)=1$ for any $F\in\rS_0(2n,\R)$.
 Since $Y\1=\0$ the linear term in the above expression in identically zero.  Hence $C$ is a critical point of $\haff_k(\cdot)$ on $\Psi_{2n}$.
 Observe next that for $X=[x_{ij}]\in \rS_0(2n,\R)$
 \[ 2\sum_{1\le i<j<p<q\le 2n} (x_{ij}x_{pq}+x_{ip}x_{jq}+x_{iq}x_{jp})=(\sum_{1\le i<j\le 2n} x_{ij})^2 -\sum_{i=1}^{2n} (\sum_{j=1}^{2n} x_{ij})^2 +\sum_{1\le i<j\le 2n} x_{ij}^2.\]
 Hence for $Y\in\Phi_{2n}$ we have the equality
 \[ \sum_{1\le i<j<p<q\le 2n} (y_{ij}y_{pq}+y_{ip}x_{jq}+y_{iq}y_{jp})=\frac{1}{2} \sum_{1\le i<j\le 2n} y_{ij}^2.\]
 So the Hessian of $\haff_k(\cdot)$ at $C$ on $\Psi_{2n}$ is strictly positive definite.

 Let $k=2$.  Since $\haff_2(C+Y)$ is  a degree $2$ polynomial the the identity
 \[\haff_2(C+Y)=\haff_2(C)+\frac{1}{2}  \sum_{1\le i<j\le 2n} y_{ij}^2 > \haff_2(C)\]
 if $Y\in \Phi_{2n}\setminus \{0\}$.  As $C+\Phi_{2n}\supset \Psi_{2n}$ we deduce that Conjecture \ref{fridcon} holds for $k=2$.
 \qed


\begin{thebibliography}{99}
 \bibitem{BS} A. Barvinok and A. Samorodnitsky, Computing the partition function for perfect matchings in a hypergraph, 
 arXiv:1009.2397v2.
 \bibitem{CCPS} W.J. Cook, W.H. Cunningham, W.R. Pulleyblank and A. Schrijver, \emph{Combinatorial Optimization}, Wiley, 1998.
 \bibitem{CPS} M. Cygan, M. Pilipczuk and R. Skrekovski, A bound on the number of perfect matchings in Klee-graphs, University of Ljubljana,
 Preprint series, vol. 47 (2009), 1105, http://www.imfm.si/preprinti/PDF/01105.pdf.
 \bibitem{Ed65} J. Edmonds, Paths, trees and flowers, \emph{Canadian Journal of Mathematics} 17 (1965), 449--467.
 \bibitem{Ed65a} J. Edmonds, Maximum mathchings and a polyhedron with $0,1$-vertices, \emph{Journal of Research
 of the National Bureau of Standards} (B) 69 (1965), 125--130.
 \bibitem{Ego}  G.P. Egorichev, Proof of the van der Waerden conjecture for permanents,
 \emph{Siberian Math. J.} 22 (1981), 854--859.

 \bibitem{ER} P. Erd\"os and A. R\'enyi, On random matrices, II,
 \emph{Studia Math. Hungar.} 3 (1968), 459-464.

 \bibitem{EKKKN} L. Esperet, F. Kardos, A. King, D. Kral and S. Norine,
 Exponentially many perfect matchings in cubic graphs, arXiv:1012.2878.

 \bibitem{Fal}  D.I. Falikman, Proof of the van der Waerden conjecture regarding the permanent of
 doubly stochastic matrix, \emph{Math. Notes Acad. Sci. USSR} 29 (1981), 475--479.

 \bibitem{Fel} W. Feller, \emph{An Introduction to Probability Theory and Its Applications}, vol. I, J. Wiley, 1958.

 \bibitem{Fr79} S. Friedland, A lower bound for the permanent of doubly stochastic matrices,
 \emph{Ann. of Math.} 110 (1979), 167-176.

 \bibitem{Fr82}  S. Friedland, A proof of a generalized van der
 Waerden conjecture on permanents,
 \emph{Lin.\ Multilin.\ Algebra} 11 (1982), 107--120.

 \bibitem{Fr11} S. Friedland, Some open problems in matchings in graphs, conference talk in \emph{Linear Algebraic Techniques in
 Combinatorics/Graph Theory}, BIRS, February 1, 2011, http://www2.math.uic.edu/$\sim$friedlan/Friedland1fFeb11.pdf

 \bibitem{FG06} S. Friedland and L. Gurvits, Generalized Friedland-Tverberg inequality: applications and extensions,
 arXiv:math/0603410v2.

 \bibitem{FG08} S. Friedland and L. Gurvits, Lower bounds for partial matchings in regular bipartite
 graphs and applications to the monomer-dimer entropy,
 \emph{Combinatorics, Probability and
 Computing}, 2008, 15pp.

 \bibitem{FKLM} S. Friedland, E. Krop,  P.H. Lundow and K. Markstr\"om,
 Validations of the Asymptotic Matching Conjectures, \emph{Journal of Statistical Physics}, 133 (2008), 513-533,
 arXiv:math/0603001v3.

 \bibitem{FKM} S. Friedland, E. Krop and K. Markstr\"om,
 On the Number of Matchings in Regular Graphs, \emph{The Electronic Journal of Combinatorics}, 15 (2008), \#R110, 1-28,
 arXiv:0801.2256v1 [math.Co] 15 Jan 2008.

 \bibitem{FL6} S. Friedland and D. Levy, A polynomial-time approximation algorithm for the number
 of $k$-matchings in bipartite graphs, 61-67,
 \emph{Mathematical papers in honour of Eduardo Marques de S\'a},
 61--67, Textos Mat. Sér. B, 39, Univ. Coimbra, Coimbra, 2006.

 \bibitem{FP}  S. Friedland and U.N. Peled, Theory of Computation of Multidimensional Entropy with an
 Application to the Monomer-Dimer Problem, \emph{Advances of Applied Math.} 34(2005), 486-522.

 \bibitem{Gur1} L. Gurvits,
 Hyperbolic polynomials approach to van der
 Waerden/Schrijver-Valiant like conjectures,
 STOC'06: Proceedings of the 38th Annual ACM Symposium on Theory of Computing,
 417--426, ACM, New York, 2006.

 \bibitem{JSV} M. Jerrum, A. Sinclair and E. Vigoda,
 A polynomial-time approximation algorithm for the permanent of a matrix
 with non-negative entries, \emph{J. ACM} 51 (2004), 671-697.

 \bibitem{LP} L. Lov\'asz anf M.D. Plummer, \emph{Matching Theory}, North-Holland, 1986.

 \bibitem{Sch}  A. Schrijver, Counting $1$-factors in regular
 bipartite graphs, \emph{J. Comb.\ Theory} B 72 (1998),
 122--135.

 \bibitem{Tve} H. Tverberg, On the permanent of bistochastic
 matrix, \emph{Math. Scand.} 12 (1963), 25-35.

 \bibitem{Val} L.G. Valiant, The complexity of computing the
 permanent, \emph{Theoretical Computer Science} 8 (1979), 189-201.

 \bibitem{Vor} M. Voorhoeve, A lower bound for the permanents of certain $(0,1)$-matrices, \emph{Neder. Akad. Wetensch.
 Indag. Math.} 41 (1979), 83-86.

 \bibitem{vdW} B.L. van der Waerden, Aufgabe 45, \emph{Jber Deutsch. Math.-Vrein.}
 35 (1926), 117.


 \end{thebibliography}
\end{document}